\documentclass{amsart}

\usepackage{fancyhdr}
\usepackage{lastpage}
\usepackage{stmaryrd,yhmath}

\newcommand{\bdis}{\begin{displaymath}}
\newcommand{\edis}{\end{displaymath}}
\newcommand{\be}{\begin{equation}}
\newcommand{\ee}{\end{equation}}
\newcommand{\mbb}{\mathbb}
\newcommand{\mcal}{\mathcal}

\newcommand{\vp}{\varphi}

\newtheorem{lemma}[]{Lemma}

\theoremstyle{definition}

\theoremstyle{remark}
\newtheorem{remark}[]{Remark}

\newtheorem*{mydef1}{{\bf Theorem}}

\numberwithin{equation}{section}

\begin{document}

\title{Jacob's ladders and the asymptotic formula for the integral of the eight order expression $|\zeta(1/2+i\vp_2(t))|^4|\zeta(1/2+it)|^4$}

\author{Jan Moser}

\address{Department of Mathematical Analysis and Numerical Mathematics, Comenius University, Mlynska Dolina M105, 842 48 Bratislava, SLOVAKIA}

\email{jan.mozer@fmph.uniba.sk}


\begin{abstract}
It is proved in this paper that there is a fine correlation between the values of $|\zeta(1/2+i\vp_2(t))|^4$ and $|\zeta(1/2+it)|^4$ where
$\vp_2(t)$ stands for the Jacob's ladder of the second order. This new asymptotic formula cannot be obtained in known theories of Balasubramanian,
Heath-Brown and Ivic.
\end{abstract}

\maketitle

\section{Results}

Let $\mu(y)\in C^\infty([y_0,\infty))$ is a monotonically increasing (to $+\infty$) function, and let it obey $\mu(y)\geq 4y\ln y$. Similarly to
\cite{3}, (3.1)-(3.9) we obtain that there exists an unique solution $x_\mu(T)=\vp_2(T;\mu)=\vp_2(T),\ T\geq T_0[\mu]$ to the
nonlinear integral equation
\be \label{1.1}
\int_0^{\mu[x(T)]}Z^4(t)e^{-\frac{t}{x(T)}}{\rm d}t=\int_0^T Z^4(t){\rm d}t.
\ee

\begin{remark}
The function $\vp_2(T)$ is to be named the \emph{Jacob's ladder of the second order}. This function obeys the following properties
\begin{itemize}
\item[(a)] it is increasing for $T\geq T_0$,
\item[(b)] if $T=\gamma$ is a zero of the function $\zeta(1/2+iT)$ of the order $n(\gamma)$ then
\bdis
\vp_2^\prime(\gamma)=\vp_2^{\prime\prime}(\gamma)=\dots=\vp_2^{(4n)}(\gamma)=0,\ \vp_2^{(4n+1)}(\gamma)\not=0 ,
\edis
\item[(c)] if
\bdis
\Phi_2(y)=\int_0^{\mu(y)}Z^4(t)e^{-\frac{t}{y}}{\rm d}t,\ y=\vp_2(T)
\edis
then
\be \label{1.2}
Z^4(t)=\Phi^\prime_2[\vp_2(T)]\frac{{\rm d}\vp_2(T)}{{\rm d}T};\ \Phi_2^\prime=\frac{{\rm d}\Phi_2}{{\rm d}\vp_2},\ T\geq T_0 ,
\ee
where
\be \label{1.3}
\Phi_2^\prime(y)=\frac{1}{y^2}\int_0^{\mu(y)}tZ^4(t)e^{-\frac{t}{y}}{\rm d}t+Z^4[\mu(y)]e^{-\frac{\mu(y)}{y}}\frac{{\rm d}\mu(y)}{{\rm d}y} .
\ee
\end{itemize}
\end{remark}

The following theorem holds true

\begin{mydef1}
If
\be \label{1.4}
[T,T+U]=\vp_2\left\{ \left[\mathring{T},\widering{T+U}\right]\right\},\ U=T^{13/14+2\epsilon} ,
\ee
then
\be \label{1.5}
\int_{\mathring{T}}^{\widering{T+U}}\left| \zeta\left(\frac{1}{2}+i\vp_2(t)\right)\right|^4\left| \zeta\left(\frac{1}{2}+it\right)\right|^4\sim
\frac{1}{4\pi^4}U\ln^8 T,\ T\to\infty .
\ee
\end{mydef1}

\begin{remark}
The formula (1.5) is the first asymptotic formula in the theory of the Riemann zeta-function for the eight order expression
\bdis
\left|\zeta\left(\frac{1}{2}+i\vp_2(t)\right)\right|^4\left|\zeta\left(\frac{1}{2}+it\right)\right|^4=Z^4[\vp_2(t)]Z^4(t) ,
\edis
where
\bdis
Z(t)=e^{i\vartheta(t)}\zeta\left(\frac{1}{2}+it\right),\ \vartheta(t)=-\frac{t}{2}\ln\pi+\text{Im}\ln\Gamma\left(\frac{1}{4}+i\frac t2\right) .
\edis
It is clear that the formula (1.5) cannot be obtained by complicated methods of Balasubramanian, Heath-Brown and Ivic (see \cite{2}).
\end{remark}

This paper is a continuation of the series of works \cite{3}-\cite{9}.

\section{Consequences of the Titchmarsh-Atkinson formula}

Titchmarsh has proved in 1928 the following formula (see \cite{10}, pp. 137, 141, \cite{11}, p. 143)
\be \label{2.1}
\int_0^\infty Z^4(t)e^{-\delta t}{\rm d}t\sim \frac{1}{2\pi^2}\frac{1}{\delta}\ln^4\frac{1}{\delta} .
\ee
Let us remind the Titchmarsh-Atkinson formula (see \cite{11}, p. 145)
\begin{eqnarray} \label{2.2}
& &
\int_0^\infty Z^4(t)e^{-\delta t}{\rm d}t=
\frac{1}{\delta}\left( A\ln^4\frac 1\delta +B\ln^3\frac 1\delta+C\ln^2\frac 1\delta +D\ln\frac 1\delta +E\right)+ \\
& &
\mcal{O}\left\{\left(\frac 1\delta \right)^{13/14+\epsilon}\right\},\ A=\frac{1}{2\pi^2} \nonumber
\end{eqnarray}
which improved the Titchmarsh formula (2.1). The following lemma is true.

\begin{lemma}
\begin{eqnarray} \label{2.3}
& &
\int_0^TZ^4(t){\rm d}t=\vp_2(T)\left\{ A\ln^4\vp_2(T)+B\ln^3\vp_2(T)+C\ln^2\vp_2(T)+\right. \\
& &
\left. +D\ln\vp_2(T) +E \right\}+\mcal{O}\left\{ \left( \vp_2(T)\right)^{13/14+\epsilon}\right\},\ A=\frac{1}{2\pi^2} . \nonumber
\end{eqnarray}
\end{lemma}
\begin{proof}
Similarly to \cite{3}, (4.3)-(4.6) we have
\begin{eqnarray*}
& &
Z^4(t)e^{-\frac{\delta}{2}t}<te^{-\frac \delta 2 t}=f_2(t;\delta)\leq f_2\left(\frac 2\delta;\delta\right)=\frac{2}{e\delta}, \\
& &
\int_U^\infty Z^4(t)e^{-\frac \delta 2 t}e^{-\frac \delta 2 t}{\rm d}t<\frac{4}{e\delta^2}e^{-\frac \delta 2 U} .
\end{eqnarray*}
The value $U=\mu(1/\delta)$ is to be chosen by the following rule
\bdis
\frac{4}{e\delta^2}e^{-\frac \delta 2 U}\leq \frac{1}{\sqrt{\delta}}\ \Rightarrow\
\mu\left(\frac 1 \delta\right)\geq\frac{4}{\delta}\ln\frac 1\delta > \frac 2\delta\ln\frac{4}{e\delta^{3/2}} .
\edis
Now (2.2) implies
\begin{eqnarray} \label{2.4}
& &
\int_0^{\mu(1/\delta)}Z^4(t)e^{-\delta t}{\rm d}t=\frac{1}{\delta}
\left( A\ln^4\frac 1\delta+B\ln^3\frac 1\delta+C\ln^2\frac 1\delta+\right. \\
& &
\left. +D\ln\frac 1\delta+E\right)+\mcal{O}\left\{ \left(\frac 1\delta\right)^{13/14+\epsilon}\right\},\
\mu\left(\frac 1\delta\right)\geq\frac 4\delta\ln\frac 1\delta . \nonumber
\end{eqnarray}
Since (see (1.1), compare \cite{3}, (3.3))
\bdis
\int_0^{\mu(1/\delta)}Z^4(t)e^{-\delta t}{\rm d}t=\int_0^{M_2(y)}Z^4(t){\rm d}t ,
\edis
and $\frac{1}{\delta}=y=\vp_2(T),\ M_2[\vp_2(T)]=T$, then (2.4) implies (2.3).
\end{proof}

\section{The asymptotic formula $\vp_2(T)\sim T$}

The following lemma is true.

\begin{lemma}
\be \label{3.1}
\vp_2(T)-T=\mcal{O}\left(\frac{T}{\ln T}\right),\ T\to\infty .
\ee
\end{lemma}
\begin{proof}
In 1924 Ingham has proved the following formula (see \cite{1}, p. 277, \cite{11}, p. 125)
\be \label{3.2}
\int_0^T Z^4(t){\rm d}t=\frac{1}{2\pi^2}T\ln^4 T+\mcal{O}(T\ln^3 T) .
\ee
Let us remind the Ingham-Heath-Brown formula (see \cite{2}, p. 129)
\be \label{3.3}
\int_0^T Z^4(t){\rm d}t=T\sum_{k=0}^4 C_k\ln^{4-k} T+\mcal{O}(T^{7/8+\epsilon}) ,
\ee
which improved the Ingham formula (3.2); $C_0=\frac{1}{2\pi^2}$ is the Ingham constant. Putting $T=M_2(y); \vp_2(T)=\vp_2[M_2(y)]=y$ into eq.
(2.3) we obtain
\begin{eqnarray} \label{3.4}
& &
\int_0^{M_2(y)}Z^4(t){\rm d}t=y\left( C_0\ln^4y+B\ln^3y+C\ln^2y+D\ln y+E\right)+\mcal{O}(y^{13/14+\epsilon}) .
\end{eqnarray}
Furthermore, putting $T=\frac{4}{5}y, \frac{5}{4}y$ into eq. (3.3) and comparing with the formula (3.4) we obtain
\bdis
\frac{4}{5}y<M_2(y)<\frac{5}{4}y\ \Rightarrow\ \frac{4}{5}\vp_2(T)<T<\frac{5}{4}\vp_2(T)\ \Rightarrow\ \frac{4}{5}T<\vp_2(T)<\frac{5}{4}T ,
\edis
i.e.
\be\label{3.5}
|\vp_2(T)-T|\leq \frac 14 T .
\ee
Now, by comparing of the formulae (2.3), (3.3) (see (3.5)) we obtain
\be \label{3.6}
C_0\left\{ \vp_2(T)\ln^4 \vp_2(T)-T\ln^4 T\right\}=\mcal{O}(T\ln^3T).
\ee
Next, from (3.6) by the Taylor formula and (3.5), we obtain
\begin{eqnarray} \label{3.7}
& &
C_0\left( \ln^4T+4\ln^3T\right)\left[\vp_2(T)-T\right]+
\mcal{O}\left\{ \frac{\ln^3 \hat{T}}{\hat{T}}\left[ \vp_2(T)-T\right]^2\right\}= \\
& &
\mcal{O}(T\ln^3T),\ \hat{T}=\mcal{O}(T) . \nonumber
\end{eqnarray}
Finally, we obtain (3.1) from (3.7).
\end{proof}

\section{Lemma about $\Phi^{\prime\prime}_{y^2}[\vp_2(T)]$}

By (1.3) we have
\be \label{4.1}
\Phi^{\prime\prime}_{y^2}(y)=J+Q,
\ee
where
\begin{eqnarray} \label{4.2}
& &
J=\frac{1}{y^3}\int_0^{\mu(y)}\left( \frac{t^2}{y}-2t\right)e^{-\frac{t}{y}}Z^4(t){\rm d}t , \\
& & \label{4.3}
Q=e^{-\frac{\mu(y)}{y}}
\left\{ \frac{2}{y^2}\mu(y)\frac{{\rm d}\mu(y)}{{\rm d}y}Z^4[\mu(y)]+4\left(\frac{{\rm d}\mu(y)}{{\rm d}y}\right)^2Z^3[\mu(y)]Z'[\mu(y)]\right.- \\
& &
\left. -\frac{1}{y}\left(\frac{{\rm d}\mu(y)}{{\rm d}y}\right)^2Z^4[\mu(y)]+\frac{{\rm d}^2\mu(y)}{{\rm d}y^2}Z^4[\mu(y)]\right\} . \nonumber
\end{eqnarray}
The following lemma is true.

\begin{lemma}
If $\mu(y)=4y\ln y$ then
\be\label{4.4}
\Phi^{\prime\prime}_{y^2}[\vp_2(T)]=\mcal{O}\left\{\frac{1}{T}\ln^4T(\ln\ln T)^2\right\} .
\ee
\end{lemma}
\begin{proof}
Let
\bdis
g_1(t)=\left(\frac{t^2}{y}-2t\right)e^{-\frac{t}{y}},\ t\in [0,\mu(y)] .
\edis
We apply the following elementary facts
\begin{eqnarray} \label{4.5}
& &
\min\{ g_1(t)\}=g_1[(2-\sqrt{2})y]=-2(\sqrt{2}-1)e^{-2+\sqrt{2}}y, \nonumber \\
& &
\max\{ g_1(t)\}=g_1[(2+\sqrt{2})y]=2(\sqrt{2}+1)e^{-2-\sqrt{2}}y , \\
& &
g_1(t)\leq g_1(y\ln\ln y)<y\frac{(\ln\ln y)^2}{\ln y},\ t\in [y\ln\ln y,4y\ln y] , \nonumber
\end{eqnarray}
and the Ingham formula (see (3.1)). First of all we have (see (4.2))
\begin{eqnarray} \label{4.6}
\frac{1}{y^3}\int_0^{y\ln\ln y}&=&\mcal{O}\left(\frac{1}{y^2}\int_0^{y\ln\ln y}Z^4(t){\rm d}t\right)=
\mcal{O}\left(\frac{1}{y^2}y\ln^4y\ln\ln y\right)= \\
&= &
\mcal{O}\left(\frac{1}{y}\ln^4 y\ln\ln y\right) , \nonumber
\end{eqnarray}
\begin{eqnarray} \label{4.7}
\frac{1}{y^3}\int_{y\ln\ln y}^{4y\ln y}&=&
\mcal{O}\left(\frac{1}{y^3}y\frac{(\ln\ln y)^2}{\ln y}\int_0^{4y\ln y}Z^4(t){\rm d}t\right)= \\
&=&
\mcal{O}\left(\frac{1}{y^2}\frac{(\ln\ln y)^2}{\ln y}y\ln^5 y\right)=
\mcal{O}\left( \frac{1}{y}\ln^4 y(\ln\ln y)^2\right) \nonumber
\end{eqnarray}
by (3.1), (4.5). Next we have (see (4.3))
\be \label{4.8}
Q(y)=\mcal{O}\left(\frac{\ln^3 y}{y^3}\right) .
\ee
Finally, we obtain (4.4) from (4.1) by (4.2), (4.3), (4.6)-(4.8).
\end{proof}

\begin{remark}
It is quite evident that our lemma (i.e. also Theorem) is true for continuous class of functions
\bdis
\mu(y)=4y^{\omega_1}\ln^{\omega_2}y , \ \omega_1,\omega_2\geq 1 .
\edis
\end{remark}

\section{Proof of the Theorem}

\subsection{}

Since
\begin{eqnarray*}
& &
G(t)=t(A\ln^4t+B\ln^3t+C\ln^2t+D\ln t+E) \ \Rightarrow \\
& &
G^\prime(t)=A\ln^4 t+\mcal{O}(\ln^3 t) ,
\end{eqnarray*}
then from (2.3) we obtain (see (3.1))
\begin{eqnarray} \label{5.1}
& &
\int_T^{T+U} Z^4(t){\rm d}t=\left[ A\ln^4\xi+\mcal{O}(\ln^3\xi)\right]\left[ \vp_2(T+U)-\vp_2(T)\right]+\mcal{O}(T^{13/14+\epsilon})  \\
& &
=\left[ C_0\ln^4\xi+\mcal{O}(\ln^3\xi)\right]U\tan[\alpha_2(T,U)]+\mcal{O}(T^{13/14+\epsilon}) , \nonumber
\end{eqnarray}
where
\begin{eqnarray} \label{5.2}
& &
U=T^{13/14+2\epsilon},\ A=C_0=\frac{1}{2\pi^2},\ \xi\in (\vp_2(T),\vp_2(T+U)) , \\
& &
\tan[\alpha_2(T,U)]=\frac{\vp_2(T+U)-\vp_2(T)}{U} , \nonumber
\end{eqnarray}
and $\alpha_2(T,U)$ is the angle of the chord of the curve $y=\vp_2(T)$ that binds the points $[T,\vp_2(T)]$ and $[T+U,\vp_2(T+U)]$. \\
Next we have (see (3.3))
\be \label{5.3}
\int_T^{T+U} Z^4(t){\rm d}t=C_0 U\ln^4 T+\mcal{O}(U\ln^3 T) .
\ee
Comparing formulae (5.1), (5.3) we obtain
\be \label{5.4}
\tan[\alpha_2(T,U)]=\frac{\ln^4\xi+\mcal{O}(\ln^3\xi)}{\ln^4 T+\mcal{O}(\ln^3 T)}+\mcal{O}\left(\frac{1}{T^\epsilon\ln^4 T}\right) .
\ee
Next we have (see (3.1), (5.2))
\begin{eqnarray*}
& &
\ln\xi=\ln\{ \vp_2(T)+\xi-\vp_2(T)\}=\ln\vp_2(T)+\mcal{O}\left\{ \frac{\vp_2(T+U)-\vp_2(T)}{\vp_2(T)}\right\}=\\
& &
=\ln T+\mcal{O}\left(\frac{1}{\ln T}\right) ,
\end{eqnarray*}
i.e.
\be \label{5.5}
\ln^4\xi=\ln^4 T+\mcal{O}(\ln^3 T);\ \ln^3\xi=\mcal{O}(\ln^3 T) .
\ee
Hence we have
\be \label{5.6}
\tan[\alpha_2(T,U)]=1+\mcal{O}\left(\frac{1}{\ln T}\right) .
\ee

\subsection{}

Next we have (see (1.2))
\begin{eqnarray*}
& &
\int_T^{T+U}Z^4(t){\rm d}t=\int_T^{T+U}\Phi_2^\prime[\vp_2(t)]{\rm d}\vp_2(t)=\\
& &
\Phi_2^\prime[\vp_2(\eta)][\vp_2(T+U)-\vp_2(T)]=\Phi_2^\prime[\vp_2(\eta)]U\tan[\alpha_2(T,U)],\ \eta\in (T,T+U),
\end{eqnarray*}
i.e. (see (5.6))
\be \label{5.7}
\int_T^{T+U}Z^4(t){\rm d}t=\Phi_2^\prime[\vp_2(\eta)]\left\{ 1+\mcal{O}\left(\frac{1}{\ln T}\right)\right\}U,\ U=T^{13/14+2\epsilon} .
\ee
Comparing formulae (5.3) and (5.7) we obtain
\be \label{5.8}
\Phi_2^\prime[\vp_2(\eta)]=C_0\ln^4 T+\mcal{O}(\ln^3 T) .
\ee
We have next (see (4.4))
\begin{eqnarray} \label{5.9}
& &
\Phi_2^\prime[\vp_2(t)]-\Phi_2^\prime[\vp_2(\eta)]=\Phi_2^{\prime\prime}[\vp_2(\rho)][\vp_2(t)-\vp_2(\eta)]= \\
& &
\mcal{O}\left\{ \frac{1}{T}\ln^4T(\ln\ln T)^2|\vp_2(t)-\vp_2(\eta)|\right\},\ t,\eta,\rho\in [T,T+U] . \nonumber
\end{eqnarray}
Since (see (3.1))
\bdis
\vp_2(t)-\vp_2(\eta)=t-\eta+\mcal{O}\left(\frac{T}{\ln T}\right)=\mcal{O}\left(\frac{T}{\ln T}\right)
\edis
then we obtain from (5.9) by (5.8)
\begin{eqnarray} \label{5.10}
& &
\Phi^\prime[\vp_2(t)]=\Phi_2^\prime[\vp_2(\eta)]+\mcal{O}\left\{ \ln^3T(\ln \ln T)^2\right\}= \\
& &
C_0\ln^4 T+\mcal{O}\left\{ \ln^3T(\ln\ln T)^2\right\} . \nonumber
\end{eqnarray}
Hence, we obtain from the formula (see (1.2))
\bdis
Z^4(t)=\Phi_2^\prime[\vp_2(t)]\frac{{\rm d}\vp_2(t)}{{\rm d}t},\ t\in [T,T+U]
\edis
by (5.10) the main formula
\be \label{5.11}
Z^4(t)=C_0\ln^4T\left\{ 1+\mcal{O}\left(\frac{(\ln\ln T)^2}{\ln T}\right)\right\}\frac{{\rm d}\vp_2(t)}{{\rm d}t},\
t\in [T,T+U],\ U=T^{13/14+2\epsilon} .
\ee

\subsection{}

Putting
\be \label{5.12}
\mcal{Z}^4(t)=\frac{Z^4(t)}{1+\mcal{O}\left\{ \frac{(\ln\ln T)^2}{\ln T}\right\}} ,
\ee
we obtain (see (5.11))
\be \label{5.13}
\mcal{Z}^4(t)=C_0\ln^4T\frac{{\rm d}\vp_2(t)}{{\rm d}t},\ t\in[T,T+U] .
\ee
Then from (5.13) by (5.3) we obtain the following $\mcal{Z}^4$-transformation
\begin{eqnarray} \label{5.14}
& &
\int_{\mathring{T}}^{\widering{T+U}}Z^4[\vp_2(t)]\mcal{Z}^4(t){\rm d}t=C_0\ln^4T\int_{\mathring{T}}^{\widering{T+U}}Z^4[\vp_2(t)]
\frac{{\rm d}\vp_2(t)}{{\rm d}t}{\rm d}t= \\
& &
C_0\ln^4T\int_T^{T+U}Z^4(t){\rm d}t=C_0^2U\ln^8T+\mcal{O}(U\ln^7T) . \nonumber
\end{eqnarray}
Since (see (5.12))
\begin{eqnarray} \label{5.15}
& &
\int_{\mathring{T}}^{\widering{T+U}}Z^4[\vp_2(t)]\mcal{Z}^4(t){\rm d}t= \\
& &
=\frac{1}{1+\mcal{O}\left\{ \frac{(\ln\ln \tau)^2}{\ln \tau}\right\}}\int_{\mathring{T}}^{\widering{T+U}}Z^4[\vp_2(t)]Z^4(t){\rm d}t,\
\tau\in (\mathring{T},\widering{T+U}),
\end{eqnarray}
and $T\to\infty\ \Rightarrow\ \tau\to\infty;\ \vp_2(\mathring{T})=T$, then the formula (1.5) follows from (5.14) by (5.15).

\thanks{I would like to thank Michal Demetrian for helping me with the electronic version of this work.}

\end{document}